\documentclass[11pt]{article}

\usepackage{amsfonts,amssymb,amsmath,amsthm,epsfig,euscript}
\usepackage{graphicx}

\newtheorem{theorem}{Theorem}

\long\def\symbolfootnote[#1]#2{\begingroup
\def\thefootnote{\fnsymbol{footnote}}\footnote[#1]{#2}\endgroup}

\newcommand{\cref}[1]{Corollary \ref{corollary:#1}}

\title{Counting alternating permutations with restricted prefix and suffix}

\author{
Ran Pan \\
\small Department of Mathematics\\[-0.8ex]
\small University of California, San Diego\\[-0.8ex]
\small La Jolla, CA 92093-0112. USA\\[-0.8ex]
\small \texttt{ran.pan.math@gmail.com}
\and
Jeffrey Remmel \\
\small Department of Mathematics\\[-0.8ex]
\small University of California, San Diego\\[-0.8ex]
\small La Jolla, CA 92093-0112. USA\\[-0.8ex]
\and
}

\date{
\small MR Subject Classifications: 05A15, 05E05
}

\begin{document}
\maketitle

\begin{abstract}
In this paper, we use Hasse diagrams and generating functions to count alternating permutations with restricted prefix and suffix of lengths 3 and 4.
In other words, for an alternating permutation $\sigma=\sigma_1\sigma_2\sigma_3\cdots\sigma_{n}\in S_{n}$,  we restrict length-3 prefixes $\sigma_1\sigma_2\sigma_3$  to follow certain patterns, such as $231$ and $132$, or follow certain restrictions such as $\sigma_2 \geq \max\{\sigma_1,\sigma_3\}+2$, similarly for prefixes of length 4.
We also study the enumeration of alternating permutations with restrictions on both prefix and suffix.

\end{abstract}

\section{Introduction}
A permutation of length $n$ is a rearrangement of $\{1,2,3,\cdots,n\}$. The set of all the permutations of length $n$ is denoted by $S_{n}$ and clearly, $S_{n}=n!$. It's well-known in combinatorics that a permutation $\sigma=\sigma_1\sigma_2\sigma_3\cdots\sigma_{n}\in S_{n}$ is called alternating if $\sigma_1<\sigma_2>\sigma_3<\sigma_4>\sigma_5<\cdots$. They are also called up-down permutations.  We denote the set of all the alternating permutations in $S_{n}$ by $A_{n}$. Désiré André showed that Euler numbers count the number of such permutations (\cite{And}), i.e., $|A_{n}|=E_{n}$, where $E_{n}$ is the $n$-th Euler number. $E_{n}$ has exponential generating function $f$ as follows
$$
f(x)=\sum_{n\geq 0}\frac{E_{n}}{n!}x^{n}=1+\frac1{1!}x+\frac1{2!}x^2+\frac2{3!}x^3+\frac5{4!}x^4+
\frac{16}{5!}x^5+\cdots=\sec(x)+\tan(x).
$$

\begin{figure}
\centering
\scalebox{0.35}
{\includegraphics{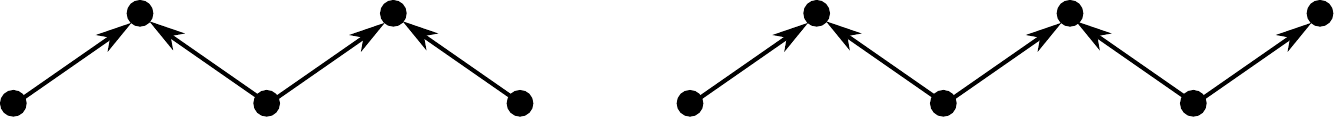}}
\caption{$A_5$ and $A_6$}
\label{fig:altperm}
\end{figure}

To better understand alternating permutations, we can use Hasse diagram to illustrate them, which is shown in Figure~\ref{fig:altperm}.
Entringer studied the enumeration of alternating permutation starting with a given value \cite{Entringer}.

In this paper, we focus on enumeration of alternating permutations with prefix or suffix that is  restricted to some certain pattern. We say a prefix or suffix has pattern $P$ if the prefix or the  suffix is order-isomorphic to $P$. For example, there is a pattern $P=1324\in A_4$ and a permutation $\sigma=24361857$. The first four elements of $\sigma$ are $2436$ which is order-isomorphic to $1324$. So we say $\sigma$ has $P$ as the prefix.

For $\sigma\in S_n$, the complement of $\sigma$ is defined as $\sigma^c=(n+1-\sigma_n,n+1-\sigma_{n-1},\cdots,n+1-\sigma_2,n+1-\sigma_1)$. For example, if $\sigma=251634$, then $\sigma^c=341625$. It it easy to see that for any pattern $P\in A_{2k}$, the number of permutations in $A_{2n}$ having $P$ as prefix is equal to the number of permutations having $P^c$ as suffix. So we will not discuss prefix and suffix separately.

This paper is mainly motivated by Remmel's~\cite{Remmel} and Bóna's~\cite{Bona}. Remmel studied pattern matching in alternating permutations and the theory in \cite{Remmel} could be extended to any minimal overlapping patterns in $A_{2n}$ and in \cite{Bona} Bóna gave results of number of minimal overlapping permutations. Part of this paper is preliminary results for a forthcoming paper about number of minimal overlapping permutations in some special families.

In Section 2, we give formulas and generating functions of the number of alternating permutations with some given restricted prefix of length 3. In Section 3, we study prefixes of length 4. In Section 4, we study alternating permutations satisfying given prefix and suffix at the same time. 
The method is based on inclusion-exclusion using Hasse diagrams and generating functions. We only consider prefix and suffix of length $3$ and $4$ and one could easily expand the discussion for longer prefix or suffix.

\section{Prefix of length $3$}
In this section, the discussion is focused on enumeration of alternating permutations with prefix of length $3$, that is, for any $\sigma=\sigma_1\sigma_2\sigma_3\sigma_4\cdots\sigma_{n+3}\in A_{n+3}$, we only take $\sigma_1\sigma_2\sigma_3$ into consideration. Prefixes that are discussed in this section include $\{\sigma_1>\sigma_3\}$, $\{\sigma_1>\sigma_3+1\}$, $\{\sigma_3>\sigma_1\}$, $\{\sigma_3>\sigma_1+1\}$, $\{\sigma_2>\max\{\sigma_1,\sigma_3\}+1\}$.

\subsection{Pattern $231$ as prefix}
\begin{figure}
\centering
\scalebox{0.35}
{\includegraphics{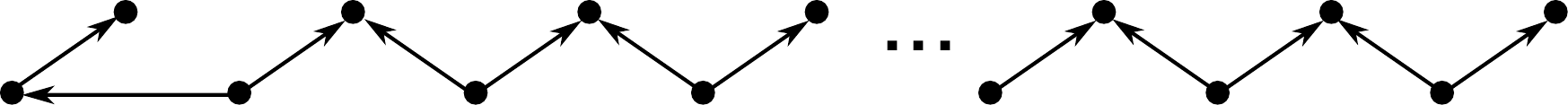}}
\caption{alternating permutations with pattern $231$ as prefix (even-length case)}
\label{fig:pre231}
\end{figure}

Suppose $\sigma=\sigma_1\sigma_2\cdots\sigma_{n+3}$ is a permutation in $A_{n+3}$ having pattern $231$ as prefix, then $\sigma_3<\sigma_1<\sigma_2$. Then we have the corresponding Hasse diagram as Figure \ref{fig:pre231} (even-length case). We denote the set of all such alternating permutations of length $n+3$ by $A^{231}(n+3)$ and clearly $|A^{{231}}(n+3)|$ is the number of linear extensions of corresponding Hasse diagram. 
We take $A^{231}(8)$ as an example and use Hasse diagrams (Figure \ref{fig:pre1324fml}) to show how inclusion-exclusion principle works here. Note that we need to be careful about choosing an appropriate edge and then apply inclusion-exclusion to it. 

\begin{figure}
\centering
\scalebox{0.3}
{\includegraphics{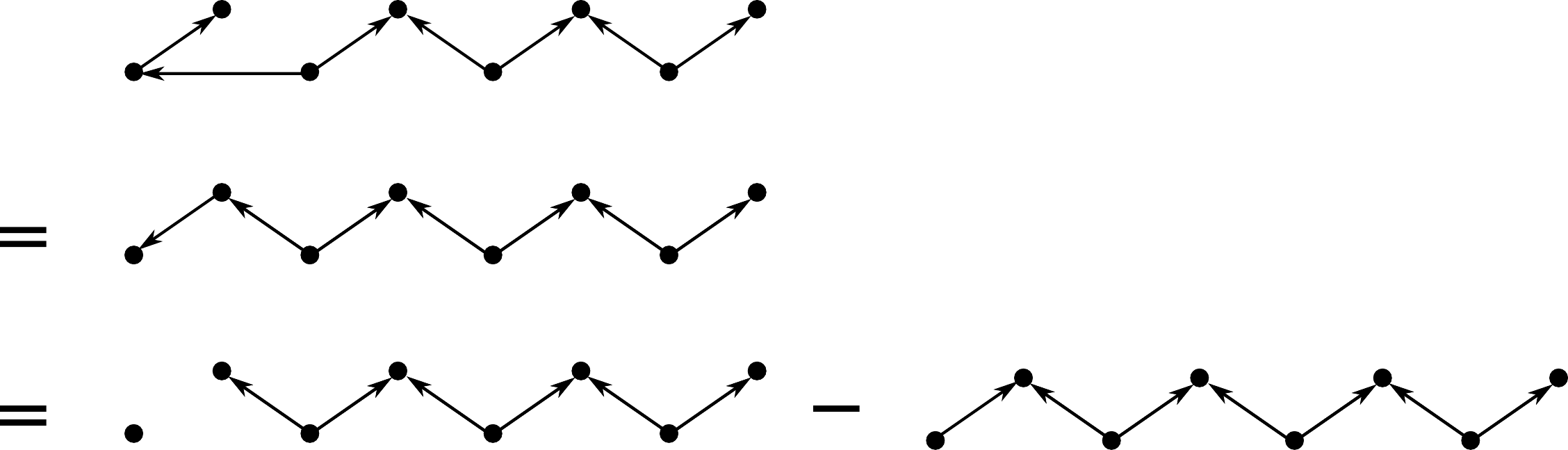}}
\caption{formula for $A^{231}(8)$}
\label{fig:pre231fml}
\end{figure}
By Figure \ref{fig:pre231fml}, since up-down permutations are essentially the same as down-up permutations, we have following formula
$$
\left|A^{231}(8)\right|=\binom{8}{1}E_7-E_{8},
$$
where $E_n$ is the $n$-th Euler number. For an arbitrary $n$, we have general formula
$$
\left|A^{231}(n+3)\right|=\binom{n+3}{1}E_{n+2}-E_{n+3}.
$$
Since $E_{n}$ has exponential generating function $f(x)=\sec x+\tan x$, the exponential generating function $f^{231}$ of $A^{{231}}(n)$ for $n\geq 3$ is
$$
f^{231}(x)=xf(x)-f(x)=(x-1)(\sec x+\tan x).
$$
A few initial terms of $|A^{231}(n)|$ for $n\geq 3$ are $1$, $3$, $9$, $35$, $155$, $791$, $4529$, $28839$, $201939$,  $\cdots$. 

In fact, observing Figure \ref{fig:pre231}, we find $f^{231}(x)$ is also the exponential generating function of up-up-down-up-down-up-down$\cdots$ permutations, which gives Sequence A034428 on OEIS \cite{Sloane} a new interpretation. 

Since pattern $231$ means $\sigma_1>\sigma_3$, a natural question could be how many alternating permutations satisfying $\sigma_1-\sigma_3\geq 2$.
It is not hard to see that number of alternating permutations of length $n$ satisfying $\sigma_1-\sigma_2=1$ is $E_{n-1}$. Therefore, for $n\geq 3$
$$
\left|A^{\sigma_1-\sigma_3\geq 2}(n)\right|=|A^{231}(n)|-E_{n-1}
$$
and the corresponding exponential generating function is
\begin{eqnarray*}
f^{\sigma_1-\sigma_3\geq 2}(x)&=&f^{231}(x)-\int_0^xf(t)dt\\
&=&(x-1)(\sec x+\tan x)+\log (1-\sin x)
\end{eqnarray*}
A few initial terms of $\left|A^{\sigma_1-\sigma_3\geq 2}(n)\right|$ for $n\geq 3$ are $0$, $1$, $4$, $19$, $94$, $519$, $3144$, $20903$, $151418$, $1188947$, $\cdots$.

\subsection{Pattern $132$ as prefix}
\begin{figure}
\centering
\scalebox{0.35}
{\includegraphics{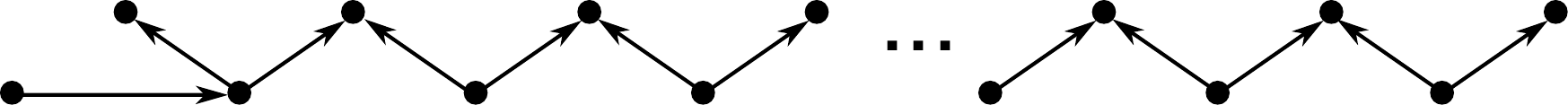}}
\caption{alternating permutations with pattern $132$ as prefix (even-length case)}
\label{fig:pre132}
\end{figure}

\begin{figure}
\centering
\scalebox{0.3}
{\includegraphics{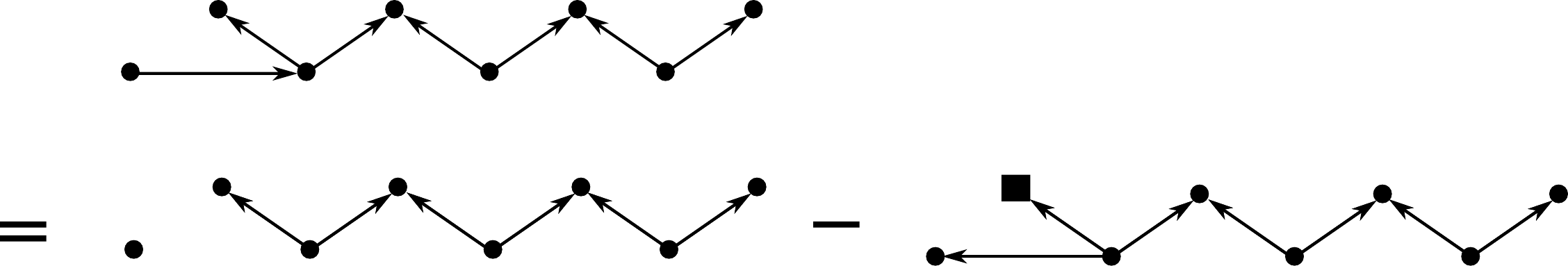}}
\caption{formula for $A^{132}(8)$}
\label{fig:pre132fml}
\end{figure}

Suppose $\sigma=\sigma_1\sigma_2\cdots\sigma_{n+3}$ is a permutation in $A_{n+3}$ having pattern $132$ as prefix, then $\sigma_1<\sigma_3<\sigma_2$. Then we have the corresponding Hasse diagram as Figure \ref{fig:pre132} (even-length case). We denote the set of all such alternating permutations of length $n+3$ by $A^{132}(n+3)$ and clearly $|A^{{132}}(n+3)|$ is the number of linear extensions of corresponding Hasse diagram. 
We take $A^{132}(8)$ as an example and use Hasse diagrams (Figure \ref{fig:pre132fml}) to show how inclusion-exclusion principle works here. Note that we need to be careful about choosing an appropriate edge and then apply inclusion-exclusion to it. 
We have two ways to put the square node in the left tail and then the new diagram we obtain is the diagram corresponding to $A^{231}(8)$. Hence, we have
$$
\left|A^{132}(8)\right|=\binom{8}{1}E_7-2\left|A^{231}(8)\right|,
$$
For an arbitrary $n$, we have general formula
$$
\left|A^{132}(n+3)\right|=\binom{n+3}{1}E_{n+2}-2\left|A^{231}(n+3)\right|.
$$
Then the exponential generating function $f^{132}$ of $A^{{132}}(n)$ for $n\geq 3$ is
\begin{eqnarray*}
f^{132}(x)&=&xf(x)-2f^{231}(x)\\
&=&(x-2(x-1))(\sec x+\tan x)\\
&=&(2-x)(\sec x+\tan x).
\end{eqnarray*}
A few initial terms of $|A^{132}(n)|$ for $n\geq 3$ are $1$, $2$, $7$, $26$, $117$, $594$, $3407$, $21682$, $\cdots$.

Actually, since $|A_n|=|A^{132}(n)|+|A^{231}(n)|$ for $n\geq 3$, there is an easier approach to get $f^{132}(x)$
\begin{eqnarray*}
f^{132}(x)&=&
f(x)-f^{231}(x)\\
&=&\left(\sec x+\tan x\right)-(x-1)\left(\sec x+\tan x\right)\\
&=&(2-x)(\sec x+\tan x).
\end{eqnarray*}
Similar to Section 2.2, we also consider number of alternating permutations satisfying $\sigma_3-\sigma_1\geq 2$.
\begin{eqnarray*}
f^{\sigma_3-\sigma_1\geq 2}(x)&=&f^{132}(x)-\int_0^xf(t)dt\\
&=&(2-x)(\sec x+\tan x)+\log (1-\sin x)
\end{eqnarray*}
A few initial terms of $\left|A^{\sigma_3-\sigma_1\geq 2}(n)\right|$ for $n\geq 3$ are $0$, $0$, $2$, $10$, $56$, $322$, $2022$, $13746$, $101332$, $806234$, $\cdots$.

\subsection{$\sigma_2>\max\{\sigma_1,\sigma_3\}+1$ in prefix}
We also explore alternating permutations $\sigma_1\sigma_2\sigma_3\cdots\sigma_n\in A_n$ such that $\sigma_2>\max\{\sigma_1,\sigma_3\}+1$. For example, we have $8$ permutations in $A_5$ having such a prefix as follows
$$
14253\phantom{123}24153\phantom{123}15243\phantom{123}15342\phantom{123}25341\phantom{123}25143\phantom{123}35142\phantom{123}35241.
$$
We split the discussion into two cases. 

Case 1. The prefix has pattern $132$. Then
$$
|A^{\sigma_2> \sigma_3+1>\sigma_1+1}(n)|=|A^{132}(n)|-|A^{\sigma_2=\sigma_3+1>\sigma_1+1}(n)|.
$$
Since $\sigma_2=\sigma_3+1$, we can combine these two nodes as one node in a corresponding Hasse diagram and then we could an up-up-down-up-down-up-down$\cdots$ permutation of length $n-1$, which implies
$$
|A^{\sigma_2> \sigma_3+1>\sigma_1+1}(n)|=|A^{132}(n)|-|A^{231}(n-1)|.
$$
Then the exponential generating function of $|A^{\sigma_2> \sigma_3+1>\sigma_1+1}(n)|$ for $n\geq 3$ is
\begin{eqnarray*}
f^{\sigma_2> \sigma_3+1>\sigma_1+1}(x)&=&f^{132}(x)-\int_0^xf^{231}(t)dt\\
&=&(2-x)(\sec x+\tan x)-\int_0^x(t-1)(\sec t+\tan t)dt.
\end{eqnarray*}
A few initial terms of $|A^{\sigma_2> \sigma_3+1>\sigma_1+1}(n)|$ for $n\geq 3$ are $0$, $1$,       $4$,           $17$,          $82$,         $439$,         $2616$,        $ 17153$, $\cdots$.

Similarly, for Case 2 where the prefix has pattern $231$.
$$
|A^{\sigma_2> \sigma_1+1>\sigma_3+1}(n)|=|A^{231}(n)|-|A^{\sigma_2=\sigma_1+1>\sigma_3+1}(n)|.
$$
Since $\sigma_2=\sigma_1+1$, we can combine these two nodes as one node in a corresponding Hasse diagram and then we could a down-up permutation of length $n-1$, which implies
$$
|A^{\sigma_2> \sigma_1+1>\sigma_3+1}(n)|=|A^{231}(n)|-E_{n-1}.
$$
Then the exponential generating function of $|A^{\sigma_2> \sigma_1+1>\sigma_3+1}(n)|$ for $n\geq 3$ is
\begin{eqnarray*}
f^{\sigma_2> \sigma_1+1>\sigma_3+1}(x)&=&f^{231}(x)-\int_0^xf(t)dt\\
&=&(x-1)(\sec x+\tan x)+\log (1-\sin x)\\
&=& f^{\sigma_1-\sigma_3\geq 2}(x)
\end{eqnarray*}
A few initial terms of $|A^{\sigma_2> \sigma_3+1>\sigma_1+1}(n)|$ for $n\geq 3$ are $0$, $1$, $4$, $19$, $94$, $519$, $3144$, $20903$, $151418$, $1188947$, $\cdots$.

Combining these two cases, we have
$$
f^{\sigma_2> \max\{\sigma_1,\sigma_3\}+1}(x)=\sec x+\tan x-\int_0^x t(\sec x+\tan t)dt
$$
and a few initial terms of $|A^{\sigma_2> \max\{\sigma_1,\sigma_3\}+1}(n)|$ for $n\geq 3$ are $0$,           $2$,           $8$,          $36$,         $176$,         $958$,        $5760$, $38056$, $\cdots$.

\section{Prefix of length $4$}
In this section, we discuss patterns of prefix of length $4$. Since $|A_4|=5$, we have $5$ permutation patterns as prefix. Fortunately, we do not need discuss them one by one because some of these permutation patterns are equivalent in the sense of enumeration.

\subsection{Pattern $1324$ as prefix}
\begin{figure}
\centering
\scalebox{0.35}
{\includegraphics{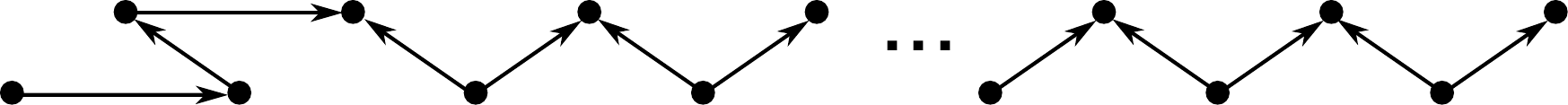}}
\caption{alternating permutations with pattern $1324$ as prefix (even-length case)}
\label{fig:pre1324}
\end{figure}

Suppose $\sigma=\sigma_1\sigma_2\cdots\sigma_{n+4}$ is a permutation in $A_{n+4}$ having $1324$ as prefix, then $\sigma_1<\sigma_3<\sigma_2<\sigma_4$. Then we have the corresponding Hasse diagram as Figure \ref{fig:pre1324} (even-length case). We denote the set of all such alternating permutations of length $n+4$ by $A^{1324}(n+4)$ and clearly $|A^{{1324}}(n+4)|$ is the number of linear extensions of corresponding Hasse diagram.

\begin{figure}
\centering
\scalebox{0.28}
{\includegraphics{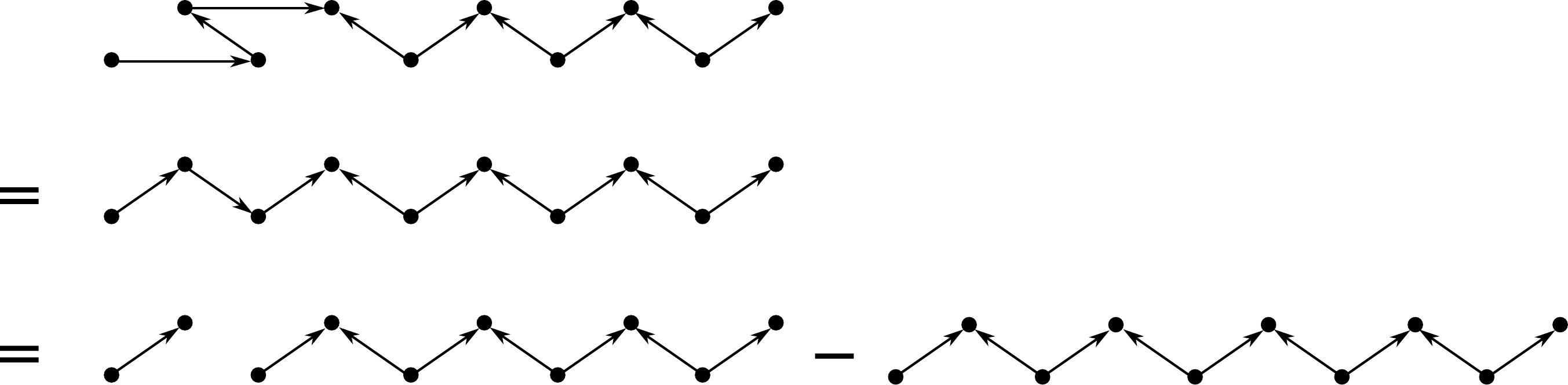}}
\caption{formula for $A^{1324}(10)$}
\label{fig:pre1324fml}
\end{figure}

We take $n=6$ as an example and use Hasse diagrams (Figure \ref{fig:pre1324fml}) to show how inclusion-exclusion principle works here. Note that we need to be careful about choosing an appropriate edge and then apply inclusion-exclusion to it. 

By Figure \ref{fig:pre1324fml}, we have following formula
$$
\left|A^{1324}(10)\right|=\binom{10}{2}E_8-E_{10},
$$
where $E_n$ is the $n$-th Euler number. For an arbitrary $n$, we have general formula
$$
\left|A^{1324}(n+4)\right|=\binom{n+4}{2}E_{n+2}-E_{n+4}.
$$
Since $E_{n}$ has exponential generating function $f(x)=\sec x+\tan x$, the exponential generating function $f^{1324}$ of $A^{{1324}}(n)$ for $n\geq 4$ is
$$
f^{{1324}}(x)=\frac12x^2f(x)-f(x)=\frac{x^2-2}2\left(\sec x+\tan x\right).
$$
A few initial terms of $|A^{1324}(n)|$ for $n\geq 4$ are $1$, $4$, $14$, $112$, $323$, $1856$,  $\cdots$.

In fact, from Figure \ref{fig:pre1324fml}, we see $|A^{{1324}}(n+4)|$ also counts the number of up-up-up-down-up-down-up-down$\cdots$ permutations. 
Actually, one worth-mentioning observation is as follows.
\begin{theorem}
Suppose $P=p_1p_2\cdots p_{k}\in A_{k}$, $Q=q_1q_2\cdots q_{k}\in A_{k}$. Let $f^{P}$ denote the exponential generating function of $|A^P(n)|$ for $n\geq k$. If $p_{k}=q_{k}$, then
$
f^{P}=f^{Q}.
$
\end{theorem}
The theorem above could be easily proved by drawing the Hasse diagrams. The theorem implies number of alternating permutations of length $2n$ with $2314$ as prefix is equal to $|A^{1324}(n)|$, that is
$$
f^{{2314}}=f^{{1324}}=\frac{x^2-2}2\left(\sec x+\tan x\right).
$$

\subsection{Pattern $1423$ as prefix}
Suppose $\sigma=\sigma_1\sigma_2\cdots\sigma_{n+4}$ is a permutation in $A_{n+4}$ having $1423$ as prefix, then $\sigma_1<\sigma_4<\sigma_2<\sigma_3$. We denote the set of all such alternating permutations of length $n+4$ by $A^{1324}(n+4)$ and the exponential generating function of $|A^{1423}(n)|$ for $n\geq 4$ by $f^{1423}(x)$. By Theorem 1, we have $f^{2413}=f^{1423}$.

We take $A^{1423}(10)$ as an example. Similarly, we apply inclusion-exclusion to the corresponding Hasse diagram. 
\begin{figure}
\centering
\scalebox{0.28}
{\includegraphics{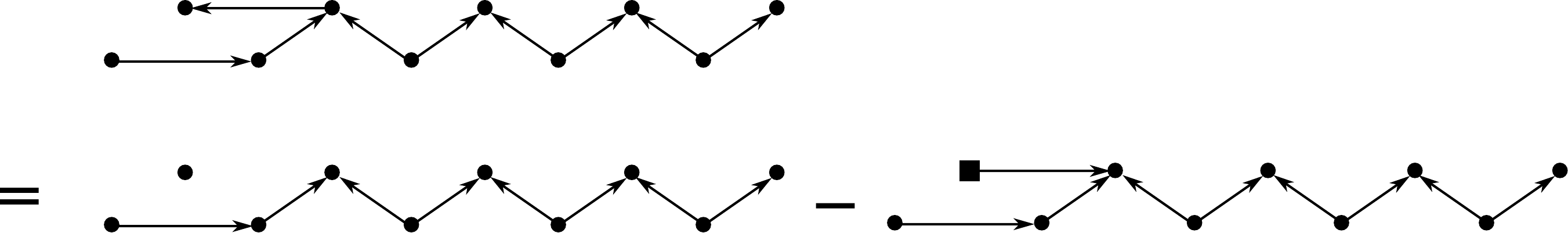}}
\caption{formula for $A^{1423}(10)$}
\label{fig:pre1423fml}
\end{figure}

For square node in the subtrahend diagram in Figure \ref{fig:pre1423fml}, we have three ways to put it and then the subtrahend becomes the Hasse diagram of $A^{1324}(10)$.The minuend diagram in in Figure \ref{fig:pre1423fml} is $A^{231}(9)$ with an isolated node. 

Then we have
$$
\left|A^{1423}(10)\right|=10|A^{231}(9)|-3\left|A^{1324}(10)\right|,
$$
and for an arbitrary $n$, we have
$$
\left|A^{1423}(n+4)\right|=(n+4)\left|A^{231}(n+3)\right|-3\left|A^{1324}(n+4)\right|.
$$
Hence the exponential generating function
\begin{eqnarray*}
f^{1423}(x)&=&xf^{231}(x)-3f^{1324}(x)\\
&=&x(x-1)\left(\sec x+\tan x\right)-3\frac{x^2-2}2\left(\sec x+\tan x\right)\\
&=&\frac{-x^2-2x+6}2\left(\sec x+\tan x\right).
\end{eqnarray*}

A few initial terms of $|A^{1423}(n)|$ for $n\geq 4$ are $1$, $3$, $12$, $53$, $271$, $1551$, $\cdots$.

Besides, by Theorem 1, we know $|A^{1423}(n)|=|A^{2413}(n)|$.

\subsection{Pattern $3412$ as prefix}
Suppose $\sigma=\sigma_1\sigma_2\cdots\sigma_{n+4}$ is a permutation in $A_{n+4}$ having $3412$ as prefix, then $\sigma_3<\sigma_4<\sigma_1<\sigma_2$. We denote the set of all such alternating permutations of length $n+4$ by $A^{3412}(n+4)$ and the exponential generating function of $|A^{3412}(n)|$ for $n\geq 4$ by $f^{3412}(x)$.

\begin{figure}
\centering
\scalebox{0.28}
{\includegraphics{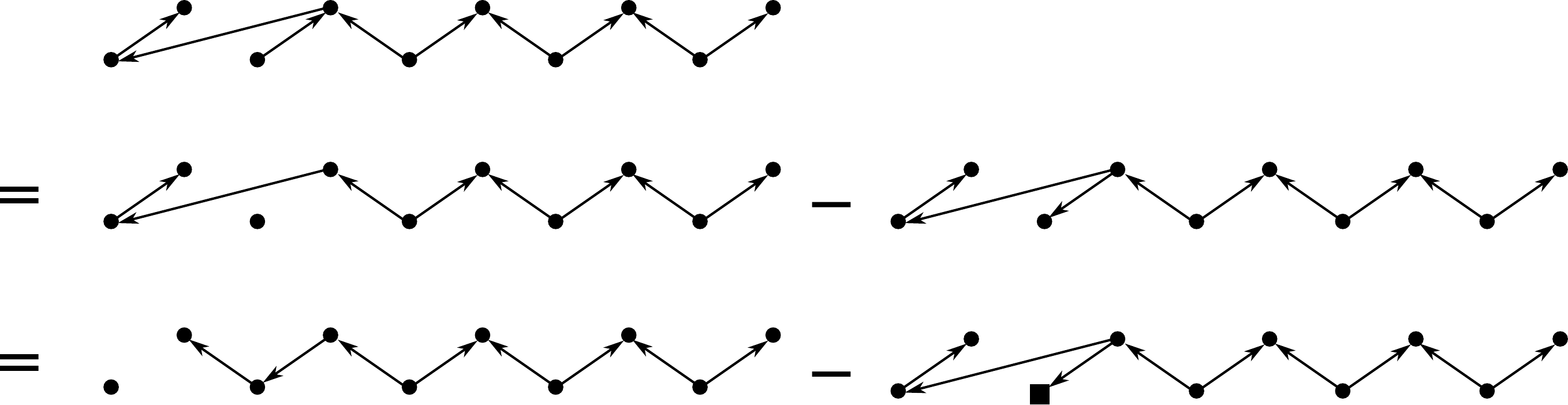}}
\caption{formula for $A^{3412}(10)$}
\label{fig:pre3412fml}
\end{figure}

\begin{figure}
\centering
\scalebox{0.28}
{\includegraphics{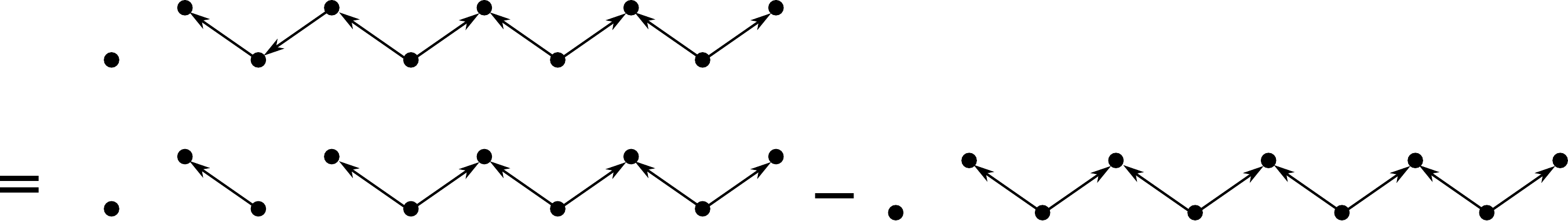}}
\caption{formula for the minuend diagram in Figure~\ref{fig:pre3412fml}}
\label{fig:pre3412min}
\end{figure}

\begin{figure}
\centering
\scalebox{0.28}
{\includegraphics{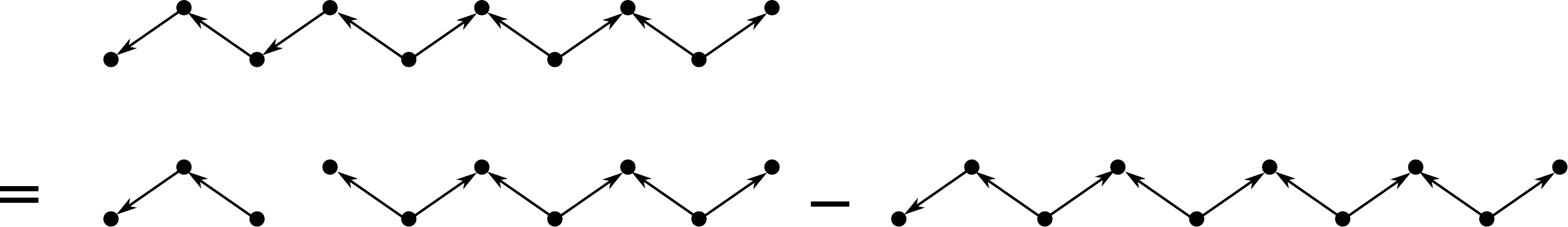}}
\caption{formula for the simplified diagram of the subtrahend diagram in Figure~\ref{fig:pre3412fml}}
\label{fig:pre3412sub}
\end{figure}

We take $A^{3412}(10)$ as an example. Similarly, as shown in Figure \ref{fig:pre3412fml} we apply inclusion-exclusion to the corresponding Hasse diagram. Then we apply inclusion-exclusion to the minuend diagram, as shown in Figure~\ref{fig:pre3412min}.  For square node in the subtrahend diagram in Figure \ref{fig:pre3412fml}, we have three ways to put it and then get a diagram shown in Figure \ref{fig:pre3412sub}. Up-down permutations and down-up permutations are essentially the same. Then combing formulas in Figure \ref{fig:pre3412min} and \ref{fig:pre3412sub}, we have
\begin{eqnarray*}
\left|A^{3412}(10)\right|&=&\left(\binom{10}1\binom{9}{2}E_7-\binom{10}{1}E_9\right)-3\left(\binom{10}3E_7-|A^{231}(10)|\right)\\
&=&2\cdot 10E_9-3E_{10}
\end{eqnarray*}
For an arbitrary $n$, we have
$$
|A^{3412}(n+4)|=2(n+4)E_{n+3}-3E_{n+4}.
$$
Hence the exponential generating function of $|A^{3412}(n)|$ for $n\geq 4$ is
\begin{eqnarray*}
f^{3412}(x)&=&2xf(x)-3f(x)\\
&=&2x\left(\sec x+\tan x\right)-3\left(\sec x+\tan x\right)\\
&=&(2x-3)\left(\sec x+\tan x\right).
\end{eqnarray*}
A few initial terms of $|A^{3412}(n)|$ for $n\geq 4$ are $1$, $2$, $9$, $38$, $197$, $1122$, $7157$, $\cdots$.

Since there are only 5 patterns in $A_4$, $f(x)=\sum_{\sigma\in A_4}f^\sigma(x)=f^{1324}(x)+f^{1423}(x)+f^{2314}(x)+f^{2413}(x)+f^{3412}(x)$ and then there is an easier way to compute $f^{3412}$
\begin{eqnarray*}
f^{3412}(x)&=&f(x)-\left(f^{1324}(x)+f^{1423}(x)+f^{2314}(x)+f^{2413}(x)\right)\\
&=&\sec x+\tan x-(x^2-2)(\sec x+\tan x)-(-x^2-2x+6)(\sec x+\tan x)\\
&=&(2x-3)(\sec x+\tan x).
\end{eqnarray*}

\subsection{Brief comparison}
This subsection mainly answers the question: which $\sigma\in A_4$ is most or least likely to be a prefix pattern, that is, we want to compare $|A^\sigma(n)|$ for $n\geq 4$ and $\sigma\in A_4$. By observation of initial terms of the sequences,  we claim that
$$
|A^{1324}(n)|=|A^{2314}(n)|>|A^{1423}(n)|=|A^{2413}(n)|>|A^{3412}(n)|.
$$
The equalities are already shown by Theorem 1. Then we only need to compare $|A^{1324}(n)|$ with $|A^{1423}(n)|$ and $|A^{1423}(n)|$ with $|A^{3412}(n)|$. From previous paragraph, we have
$$
|A^{1324}(n)|-|A^{1423}(n)|=n(n-1)E_{n-2}+nE_{n-1}-4E_n,
$$
and
$$
|A^{1423}(n)|-|A^{3412}(n)|=6E_n-\frac12n(n-1)E_{n-2}+3nE_{n-1}
$$
Since we know for Euler numbers $\lim_{n\rightarrow \infty}\frac{E_n}{E_{n-1}}=\frac{2n}{\pi}$,
\begin{eqnarray*}
\lim_{n\rightarrow \infty}|A^{1324}(n)|-|A^{1423}(n)|&=&\lim_{n\rightarrow \infty}n(n-1)E_{n-2}+nE_{n-1}-4E_n\\
&=&\frac{\pi^2+2\pi-16}4E_n\\
&\approx& 0.038 E_n
\end{eqnarray*}
and
\begin{eqnarray*}
\lim_{n\rightarrow \infty}|A^{1423}(n)|-|A^{3412}(n)|&=&\lim_{n\rightarrow \infty}6E_n-\frac12n(n-1)E_{n-2}+3nE_{n-1}\\
&=&\frac{48-\pi^2-12\pi}8E_n\\
&\approx& 0.054 E_n.
\end{eqnarray*}
Hence the following inequality
$$
|A^{1324}(n)|>|A^{1423}(n)|>|A^{3412}(n)|
$$
holds asymptotically.

\section{Prefix and Suffix of length $4$}
From this section, we will consider both the prefix and suffix. If a alternating permutation has permutation pattern of length $4$ as both its prefix and suffix, the permutation is of even length. There are $25$ combinations of prefix and suffix and fortunately some of them are equivalent. 

\subsection{Pattern $1324$ as prefix and $1324$ as suffix}
We denote the set of all  alternating permutations having $1324$ as both prefix and suffix of length $2n$ by $A^{1324}_{1324}(2n)$ and clearly, $n \geq 4$. The exponential generating function of $|A^{1324}_{1324}(2n)|$ is denoted by $f^{1324}_{1324}(x)$. 

\begin{figure}
\centering
\scalebox{0.28}
{\includegraphics{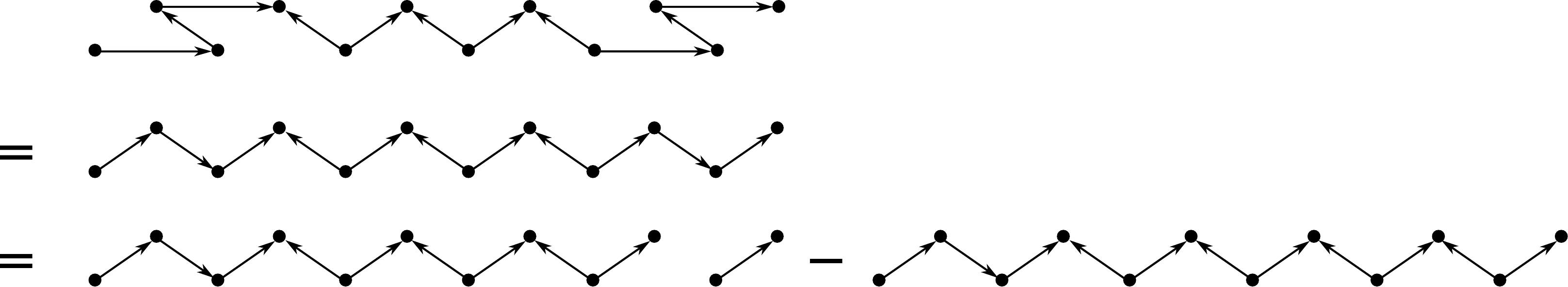}}
\caption{formula for $A^{1324}_{1324}(12)$}
\label{fig:ps1324fml}
\end{figure}

We take $A^{1324}_{1324}(12)$ as an example and the corresponding Hasse diagram is shown in Figure \ref{fig:ps1324fml} and apply inclusion-exclusion to it.

Then we have
$$
\left|A^{1324}_{1324}(12)\right|=\binom{12}2\left|A^{1324}(10)\right|-\left|A^{1324}(12)\right|
$$
and an arbitrary $n$,
$$
\left|A^{1324}_{1324}(2n+8)\right|=\binom{2n+8}2\left|A^{1324}(2n+6)\right|-\left|A^{1324}(2n+8)\right|.
$$
Note that the generating function of $|A^{1324}(2n)|$ is $\frac{x^2-2}2\sec x$. Then we obtain the exponential generating function of $|A^{1324}_{1324}(2n)|$ for $n\geq 4$,
\begin{eqnarray*}
f^{1324}_{1324}(x)&=&\left(\frac{1}2x^2\left(\frac{x^2-2}2\sec x\right)\right)-\left(\frac{x^2-2}2\sec x\right)\\
&=&\left(\frac{x^2-2}2\right)^2\sec x.
\end{eqnarray*}
A few initial terms of $|A^{1324}_{1324}(2n)|$ for $n\geq 4$ are $69$, $2731$, $147443$, $10886877$, $1059070505$, $\cdots$.

Based on Theorem 1 and considering the complement, we have
$$
f^{1324}_{1324}(x)=f^{1324}_{1423}(x)=f^{2314}_{1324}(x)=f^{2314}_{1423}(x).
$$

\subsection{Pattern $1423$ as prefix and $1423$ as suffix}
We denote the set of all  alternating permutations having $1324$ as both prefix and suffix of length $2n$ by $A^{1423}_{1423}(2n)$ and clearly, $n \geq 4$. The exponential generating function of $|A^{1423}_{1423}(2n)|$ is denoted by $f^{1423}_{1423}(x)$. 

\begin{figure}
\centering
\scalebox{0.26}
{\includegraphics{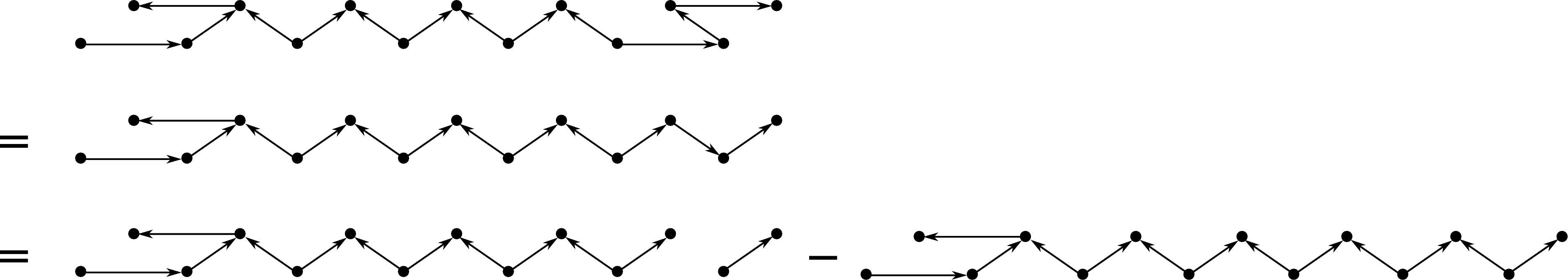}}
\caption{formula for $A^{1423}_{1423}(14)$}
\label{fig:ps1423fml}
\end{figure}

Similarly, we take $A^{1423}_{1423}(14)$ as an example and the corresponding Hasse diagram is shown in Figure \ref{fig:ps1423fml} and apply inclusion-exclusion to it.

Then we have
$$
\left|A^{1423}_{1423}(14)\right|=\binom{14}2\left|A^{1423}(12)\right|-\left|A^{1423}(14)\right|
$$
and for an arbitrary $n$,
$$
\left|A^{1423}_{1423}(2n+8)\right|=\binom{2n+8}2\left|A^{1423}(2n+6)\right|-\left|A^{1423}(2n+8)\right|.
$$
Note that the generating function of $|A^{1423}(2n)|$ is $\frac{-x^2\sec x-2x\tan x+6\sec x}2$. Then we obtain the exponential generating function of $|A^{1423}_{1423}(2n)|$ for $n\geq 4$,
\begin{eqnarray*}
f^{1423}_{1423}(x)&=&\left(\frac{1}2x^2\frac{-x^2\sec x-2x\tan x+6\sec x}2\right)-\frac{-x^2\sec x-2x\tan x+6\sec x}2\\
&=&-\frac{(x^2-2)(x^2+2x\sin x-6)\sec x}4.
\end{eqnarray*}
A few initial terms of $|A^{1423}_{1423}(2n)|$ for $n\geq 4$ are $65$, $2317$, $123543$, $9109111$, $885987557$, $\cdots$.

Based on Theorem 1 and considering the complement, we have
\begin{eqnarray*}
&&f^{1423}_{1423}(x)=f^{1423}_{1324}(x)=f^{2413}_{1423}(x)=f^{2413}_{1324}(x)\\
&=&f^{2314}_{2314}(x)=f^{1324}_{2314}(x)=f^{1324}_{2413}(x)=f^{2314}_{2413}(x).
\end{eqnarray*}

\subsection{Pattern $2413$ as prefix and $2413$ as suffix}
We denote the set of all  alternating permutations having $1324$ as both prefix and suffix of length $2n$ by $A^{2413}_{2413}(2n)$ and clearly, $n \geq 4$. The exponential generating function of $|A^{2413}_{2413}(2n)|$ is denoted by $f^{2413}_{2413}(x)$. 

\begin{figure}
\centering
\scalebox{0.28}
{\includegraphics{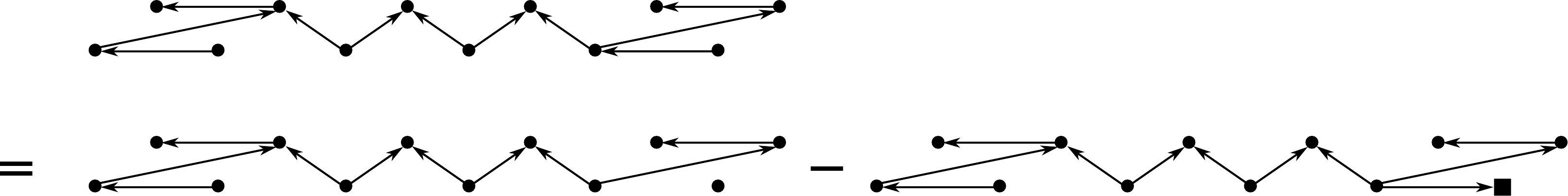}}
\caption{formula for $A^{2413}_{2413}(12)$}
\label{fig:ps2413fml}
\end{figure}

Similarly, we take $A^{2413}_{2413}(12)$ as an example and the corresponding Hasse diagram is shown in Figure \ref{fig:ps2413fml} and apply inclusion-exclusion to it.

For the minuend in Figure \ref{fig:ps2413fml}, we apply inclusion-exclusion and then find the diagram can be expressed by a combination of $A^{2413}$ and singletons. For the subtrahend in Figure \ref{fig:ps2413fml}, there are three ways to put the square node in the right tail and then $A^{2413}_{1324}(12)$ occurs. Therefore, we have the following formula
$$
\left|A^{2413}_{2413}(12)\right|=\binom{12}1\binom{11}1\left|A^{2413}(10)\right|-
\binom{12}1\left|A^{2413}(11)\right|
-3\left|A^{2413}_{1324}(12)\right|
$$
and for an arbitrary $n$,
\begin{eqnarray*}
\left|A^{2413}_{2413}(2n+8)\right|
&=&(2n+8)(2n+7)\left|A^{2413}(2n+6)\right|-
(2n+8)\left|A^{2413}(2n+7)\right|
-3\left|A^{2413}_{1324}(2n+8)\right|\\
&=&
(2n+8)(2n+7)\left|A^{1423}(2n+6)\right|-
(2n+8)\left|A^{1423}(2n+7)\right|
-3\left|A^{1423}_{1423}(2n+8)\right|.
\end{eqnarray*}
Note that the generating function of $|A^{1423}(2n)|$ is $\frac{-x^2\sec x-2x\tan x+6\sec x}2$ and the generating function of $|A^{1423}(2n+1)|$ is $\frac{-x^2\tan x-2x\sec x+6\tan x}2$. Then we have the exponential generating function of $|A^{2413}_{2413}(2n)|$ for $n\geq 4$ as follows
\begin{eqnarray*}
f^{2413}_{2413}(x)
&=&x^2\frac{-x^2\sec x-2x\tan x+6\sec x}2-x\frac{-x^2\tan x-2x\sec x+6\tan x}2-3f^{1423}_{1423}(x)\\
&=&-\frac{x\left(x\left(x^2-8\right)+\left(x^2+6\right)\sin x\right)\sec x}{2}+\frac{3(x^2-2)(x^2+2x\sin x-6)\sec x}4\\
&=&\frac{(x^2-4)\left(x^2+2x\sin x-6\right)\sec x}{4}.
\end{eqnarray*}
A few initial terms of $|A^{2413}_{2413}(2n)|$ for $n\geq 4$ are $53$, $1929$, $103287$, $7619719$, $741197433$, $\cdots$.

Based on Theorem 1 and considering the complement, we have
\begin{eqnarray*}
f^{2413}_{2413}(x)=f^{1423}_{2413}(x)=f^{1423}_{2314}(x)=f^{2413}_{2314}(x).
\end{eqnarray*}

\subsection{Pattern $3412$ as prefix and $1324$ as suffix}
 We denote the set of all such alternating permutations of length $2n$ by $A^{3412}_{1324}(2n)$ and clearly, $n \geq 4$. The exponential generating function of $|A^{3412}_{1324}(2n)|$ is denoted by $f^{3412}_{1324}(x)$.

\begin{figure}
\centering
\scalebox{0.28}
{\includegraphics{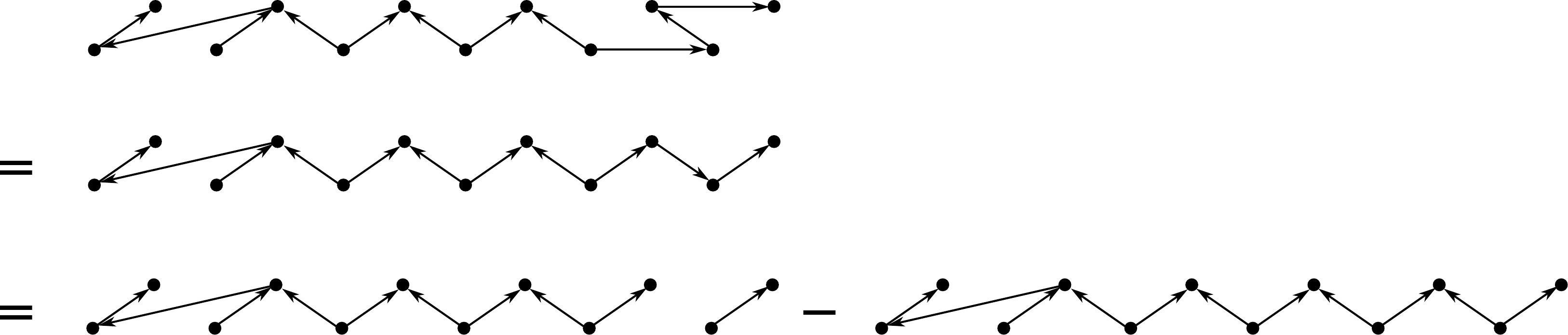}}
\caption{formula for $A^{3412}_{1324}(12)$}
\label{fig:ps3412fml}
\end{figure}

Similarly, we take $A^{3412}_{1324}(12)$ as an example and the corresponding Hasse diagram is shown in Figure \ref{fig:ps3412fml} and apply inclusion-exclusion to it. 

Then we have
$$
\left|A^{3412}_{1324}(12)\right|=\binom{12}2\left|A^{3412}(10)\right|-
\left|A^{3412}(12)\right|
$$
and for an arbitrary $n$,
$$
\left|A^{3412}_{1324}(2n+8)\right|=\binom{2n+8}2\left|A^{3412}(2n+6)\right|-
\left|A^{3412}(2n+8)\right|.
$$
Note that the generating function of $|A^{3412}(2n)|$ is $2x\tan x-3\sec x$ and hence the exponential generating function of $|A^{3412}_{1324}(2n)|$ for $n\geq 4$
\begin{eqnarray*}
f^{3412}_{1324}(x)&=&\frac12x^2\left(2x\tan x-3\sec x\right)-\left(2x\tan x-3\sec x\right)\\
&=&\left(\frac12x^2-1\right)\left(2x\tan x-3\sec x\right)\\
&=&\frac{(x^2-2)(2x\sin x-3)\sec x}{2}
\end{eqnarray*}
A few initial terms of $|A^{3412}_{1324}(2n)|$ for $n\geq 4$ are $55$, $1708$, $89649$, $6598658$, $641689451$, $\cdots$.

Based on Theorem 1 and considering the complement, we have
\begin{eqnarray*}
f^{3412}_{1324}(x)=f^{3412}_{1423}(x)=f^{1324}_{3412}(x)=f^{2314}_{3412}(x).
\end{eqnarray*}

\subsection{Pattern $3412$ as prefix and $2413$ as suffix}
We denote the set of all such alternating permutations of length $2n$ by $A^{3412}_{2413}(2n)$ and clearly, $n \geq 4$. The exponential generating function of $|A^{3412}_{2413}(2n)|$ is denoted by $f^{3412}_{2413}(x)$.

\begin{figure}
\centering
\scalebox{0.28}
{\includegraphics{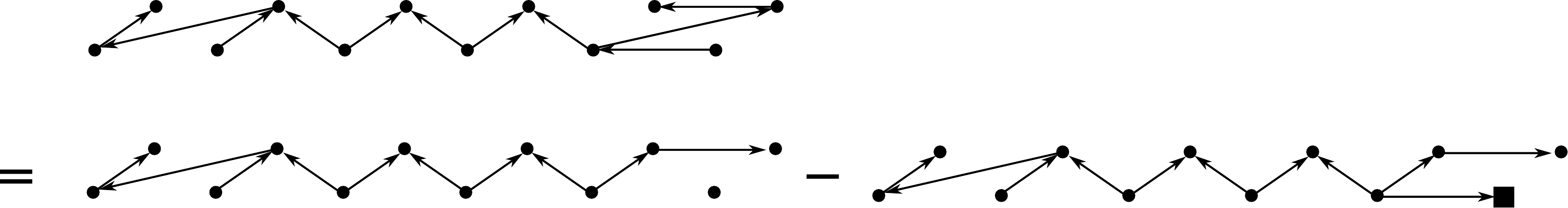}}
\caption{formula for $A^{3412}_{2413}(12)$}
\label{fig:ps34122413fml}
\end{figure}

Similarly, we take $A^{3412}_{2413}(12)$ as an example and the corresponding Hasse diagram is shown in Figure \ref{fig:ps34122413fml} and apply inclusion-exclusion to it. 

Note that there are three ways to put the square node in the right tail and then we could get a Hasse diagram corresponding to $A^{3412}_{1324}(12)$. We also need to apply inclusion-exclusion to the minuend diagram. 

Then we have formula as follows
$$
|A^{3412}_{2413}(2n+8)|=(2n+8)(2n+7)|A^{3412}(2n+6)|-(2n+8)|A^{3412}(2n+7)|-3|A^{3412}_{1324}(2n+8)|
$$
Hence the corresponding exponential generating function is as follows,
$$
f^{3412}_{2413}(x)=-\frac{\left(x^2+2(x^2-9)x\sin x+18\right)\sec x}2
$$
A few initial terms of $|A^{3412}_{2413}(2n)|$ for $n\geq 4$ are $35$, $1386$, $74745$, $5518084$, $536785743$, $\cdots$.

Based on Theorem 1 and considering the complement, we have
\begin{eqnarray*}
f^{3412}_{2413}(x)=f^{3412}_{2314}(x)=f^{1423}_{3412}(x)=f^{2413}_{3412}(x).
\end{eqnarray*}

\subsection{Pattern $3412$ as prefix and $3412$ as suffix}

Since 
$$
\sum_{\sigma,\mu\in A_4}f^\sigma_\mu(x)=\sec x,
$$
we have
\begin{eqnarray*}
f^{3412}_{3412}(x)&=&\sec x-4f^{1324}_{1324}(x)-8f^{1423}_{1423}(x)-4f^{2413}_{2413}(x)-4f^{3412}_{1324}(x)\\
&=&(4x^2-12x\sin x+9)\sec x
\end{eqnarray*}
A few initial terms of $|A^{3412}_{3412}(2n)|$ for $n\geq 4$ are $17$, $969$, $53925$, $3994741$, $388747161$, $\cdots$.

So far all the $25$ combinations of prefix and suffix of length $4$ have been discussed.

\section{Further discussion and special note}
Other than a single pattern, we could also restrict the prefix or suffix to a set of patterns. For example, $S=\{\sigma=\sigma_1\sigma_2\sigma_3\sigma_4\sigma_5\in A_5: \sigma_1<\sigma_5<\sigma_2\}$, we could ask how many alternating permutations that have $S$ as prefix.

In this paper, we only take patterns of length $4$ into account. But actually we could also consider  prefix and suffix of different lengths which could be even or odd. The idea to solve it is also based inclusion-exclusion. But more steps of inclusion-exclusion may be needed for a longer prefix or suffix.

More general problems could be study of statistics, such as descents, inversions in permutations  with given prefix or suffix.

The majority of this work was done by 2017 before the second author suddenly passed away. To honor and remember the great career of Remmel, the first author recently continued to work on the draft and decided to make it public to more readers.

\end{document}